\documentclass[a4paper,11pt]{article}
\usepackage{amsmath,amsfonts,amssymb,amscd,amsthm,epsfig}
\usepackage[mathscr]{euscript}
\usepackage{color}
\newtheorem{thm}{Theorem}[section]
\newtheorem*{*thm}{Theorem}
\newtheorem{lemma}[thm]{Lemma}
\newtheorem{prop}[thm]{Proposition}
\newtheorem{corr}[thm]{Corollary}
\theoremstyle{definition}
\newtheorem*{nota}{Notation and terminology}
\newtheorem{dfn}[thm]{Definition}
\newtheorem{exmple}[thm]{Example}
\newtheorem{exmples}[thm]{Examples}

\theoremstyle{remark}

\newtheorem*{rmq}{\textit{Remark}}
\newtheorem{rmk}[thm]{\textit{Remark}}
\renewcommand{\proof}{\noindent\textit{Proof}\/: \,\,}
\newcommand{\C}{{\mathbb{C}}}
\newcommand{\bD}{{\mathbb{D}}}
\newcommand{\Q}{{\mathbb{Q}}}
\newcommand{\R}{{\mathbb{R}}}

\newcommand{\Z}{{\mathbb{Z}}}

\newcommand{\bP}{{\mathbb{P}}}
\newcommand{\comp}{\raise1pt\hbox{{$\scriptscriptstyle\circ$}}}
 
\newcommand{\underint}{\int_{\raise-4pt\hbox{\hskip-8pt $-$}}}
\newcommand{\overint}{\int^{\raise3pt\hbox{\hskip-7pt $-$}
}\hskip -4pt}

\def\vc#1{\mathbf{#1}}

\def\lset{\{}  
\def\rset{\}}  
\def\set#1{\lset#1\rset} 
\def\st{\mid}   
\def\sett#1#2{\lset #1 \st #2 \rset}

\newcommand\Tr{{}^{\mathsf{T}}\kern-0.9pt} 
\def\id{\mathop{\rm id}\nolimits}
\newcounter{lijstc}
{\end{list}}
{\end{list}}

\def\arrow(#1,#2)\dir(#3,#4)\long#5{\put(#1,#2){\vector(#3,#4){#5}}}

\newcommand\Grid{\setbox13=\vbox to 5\unitlength{\hrule width 109mm\vfill}
\setbox13=\vbox to 65mm{\offinterlineskip\leaders\copy13\vfill\kern-1pt\hrule}
\setbox14=\hbox to 5\unitlength{\vrule height 65mm\hfill}
\setbox14=\hbox to 109mm{\leaders\copy14\hfill\kern-2mm\vrule height 65mm}
\ht14=0pt\dp14=0pt\wd14=0pt \setbox13=\vbox to
0pt{\vss\box13\offinterlineskip\box14} \wd13=0pt\box13}

\def\mapright#1{\mathop{\vbox{\ialign{
                                ##\crcr
    ${\scriptstyle\hfil\;\;#1\;\;\hfil}$\crcr
 \noalign{\kern2pt\nointerlineskip}
    \rightarrowfill\crcr}}\;}}

\def\mapleft#1{\mathop{\vbox{\ialign{
                                ##\crcr
    ${\scriptstyle\hfil\;\;#1\;\;\hfil}$\crcr
 \noalign{\kern2pt\nointerlineskip}
    \leftarrowfill\crcr}}\;}}

\newcommand\rarrow[3]{\smash{\mathop{\hbox to#3{\rightarrowfill}}\limits
^{\scriptstyle#1}_{\scriptstyle#2}}}
 
\newcommand\larrow[3]{\smash{\mathop{\hbox to#3{\leftarrowfill}}\limits
^{\scriptstyle#1}_{\scriptstyle#2}}}

\def\into{\hookrightarrow}

\newcommand\so[1]{\operatorname{SO}({#1})}

\newcommand\ogr[1]{\operatorname{O}({#1})}

\newcommand\aut[1]{\operatorname{Aut}({#1})}

\newcommand\disc[1]{\mathop{\rm discr}\nolimits(#1)}

\def\OO{{\mathcal O}}
\def\upvect{{\scriptscriptstyle\uparrow}}
\def\ogrup#1{\operatorname{O}^\upvect({#1})}

\def\ogrp#1{\operatorname{O}^+({#1})}
\def\sogrp#1{\operatorname{SO}^+({#1})}
 \def\extgrup#1{\operatorname{\tilde O}^\upvect({#1})}
\def\half{\frac 1 2} 
\begin{document}
\title{Automorphs of indefinite  binary  quadratic forms  and   K3-surfaces with Picard number $2$}
\author{Federica GALLUZZI,Giuseppe LOMBARDO and Chris PETERS}
\maketitle

\abstract{Every  indefinite binary form occurs as the Picard lattice of some K3-surface. The group of its isometries, or automorphs, coincides  with the automorphism group of the K3-surface, but only up to finite groups. The classical theory of automorphs for binary forms can then be applied to study these automorphism groups. Some extensions  of  the classical theory are needed  to  single out the orthochronous automorphs, i.e. those that conserve  the ``light cone''. Secondly, one needs to study in detail the effect of the automorphs on the discriminant group.  The result is a precise description of  all possible automorphism groups of  ``general'' K3's with Picard number two.}

\section*{Introduction} 

A K3-surface is a simply connected projective surface with trivial canonical bundle. Despite this abstract definition, K3's have been classified in detail. See for instance \cite[Chap.~8]{BPV} and the literature cited there. In \S~\ref{sec3} we collect the necessary material.

The general theory of automorphisms of K3-surfaces is largely due to Nikulin, cf.  \cite{N1,N2,N3}. The case of Picard number $1$ turns out to be quite easy to deal with. The automorphism group is finite and almost always the identity \cite{N3}.  The question of which finite groups are possible has been answered in detail by Nikulin and this question has been further investigated by Mukai \cite{Muk} and Kond\=o \cite{Kon1,Kon2}. Further attention has also been given to those groups that act trivially on the transcendental lattice, the \emph{symplectic} automorphisms, \cite{GS,GS2}.

The case of Picard number $2$ has been briefly touched upon in \cite{PSS}  where Severi's example \cite{Sev} of  a K3 with the infinite dihedral group 
as automorphism group 
is put into perspective by tying it in with the classical theory of binary quadratic forms.  Furthermore, some special cases of K3's with Picard group of rank $2$ appear in the literature \cite{W, Bi, Ge, GL}.

This note can be viewed as a systematic study of possible automorphism groups of K3's with Picard number $2$. The above examples will be placed within this general context. 

As one can surmise, number theory of integral bilinear forms 
{plays a central role through the solutions for the classical Pell equations}. Since many algebraic geometers are not very familiar with this classical theory we have  collected some of these results in \S~\ref{sec2}.  They are complemented with  new results, e.g in   \S~\ref {ssec:ambi}--\ref{repsnotzero}.
In  \S~\ref{sec3} these are  applied to  automorphism groups of K3-surfaces. It turns out that for  ``most'' K3-surfaces with Picard number $2$ the  automorphism group can be determined  in terms of  explicit number theoretic properties of the intersection form on the Picard group. 
In \S~\ref{sec4} we  give some examples for which the number theoretic data can be made explicit and we explain where the examples previously treated in the  literature   fit in.

The first two authors would like to thank Bert van Geemen for several fruitful discussions.

\begin{nota} A quadratic  form  \[ q(x_1,\dots,x_n)=\sum_{i\le j} q_{ij}x_ix_j\] in $n$ variables with coefficients in a field $k$ of characteristic $\not=2$  determines and is  determined by  the  bilinear form  $Q$  obtained from it  by polarization.  By convention  $Q$  is obtained from $q$ by placing   $q_{ii}$ on the $i$-th diagonal entry  and $\half q_{ij}$ on the $(i,j)$-th and $(j,i)$-th entry. Its determinant  $d(q)$ is called the \emph{discriminant}  of $q$.  It is well defined up to squares in $k$.  If $d(q)=  \pm 1$ the form   is called \emph{unimodular}. A word of warning: conventionally, for a quadratic form $ax^2+bxy+cy^2$ in two variables, its discriminant $b^2-4ac$ is the negative of the discriminant of the associated bilinear form!

Note that if  a quadratic form  has integral coefficients,  its associated bilinear form   has half-integral off-diagonal elements.  Nevertheless such quadratic forms are called  \emph{integral}. 

If $X$ is any complex projective variety we let $\aut X$ be its group of biholomorphic automorphisms.
A group $G$ generated by $a,b,c,\dots$ is denoted $\langle a,b,c,\dots\rangle$.
The infinite dihedral group $\bD_\infty=\Z/2\Z \ast \Z/2\Z=\langle  s, t\rangle $ is the group generated by two non-commuting involutions $s$ and $t$.
\end{nota}

\section{Automorphism groups  of  lattices} \label{sec1}

A \emph{lattice} is a pair $(S,Q)$ of a free finite rank $\Z$-module $S$ together with a bilinear form $Q:S\times S\to \Z$.   So  $S=\Z^r$  with   the standard basis yields a bilinear form   with integral coefficients.  A lattice is \emph{even}  if  $Q(\mathbf{v},\mathbf{v})\in 2\Z$ 
for all $\mathbf{v}\in S$.  %
If $S=\Z^r$ this  means that the diagonal entries of $Q$ with respect to the standard basis are even. A lattice which is not even is called \emph{odd}.   For a bilinear form $Q$ and $m\in \Z$, the form $mQ$ denotes $(\mathbf{v},\mathbf{w})\mapsto mQ(\mathbf{v},\mathbf{w})$. For instance, if  $q$ is integral, the form $2Q$  determines an even lattice and conversely.
  
 An even indefinite unimodular form is uniquely determined by its parity (i.e. if it is even or not) and its signature.   
 See e.g. \cite{Se}. 
 For instance, the  hyperbolic plane $H$   given by 
 \[
 ( \mathbf{u},\mathbf{v})  \mapsto \Tr \mathbf{u} \begin{pmatrix} 0 & 1 \\ 1& 0 \end{pmatrix} \mathbf{v}.
\]
is   the unique unimodular even lattice  of signature $(1,1)$ which is denoted by the same symbol.  
The group of isometries of a lattice $(S,Q)$ is denoted  either as  $\ogr S$ or as $\ogr  Q$. Those of determinant $1$ are denoted $\so S$ or $\so Q$.

For any lattice $(S,Q)$ the dual lattice $S^*$  is defined by 
\[
S^*= \sett{x\in S\otimes \Q}{Q(x,y)\in \Z \text{ for all } y\in 
{S}}.
\]
We have  $S\subset S^*$ and the quotient
\[
\disc S =S^*/S \quad \text{ the \emph{discriminant group} of $S$} 
\]
 is a finite abelian group of order equal to the absolute value of the discriminant $|d(Q)|$.  An isometry $g$ of $S$ induces a group automorphism on $S^*/S$ which will be noted $\bar g$.

Now we assume that  $(S,Q)$ is a  lattice of signature $(1,k)$ with $k\ge 1$.  In this situation the \emph{light cone}
\begin{equation}\label{eqn:cone}
 \sett{ x\in S\otimes_\Z\R}{ Q(x,x)>0}= C \coprod   -C 
 \end{equation}
decomposes in two connected components $C$ and $-C$. We introduce the \emph{orthochronous Lorentz group}  and the  \emph{special orthochronous Lorentz group}
\[
 \ogrp{S}= \sett{g\in \ogr{S}}{g(C )=C} ,\quad  \sogrp{S}= \ogrp{S}\cap \so{S}
 .
\]
A \emph{root} of $S$ is a vector $\mathbf{d}$ with $Q(\mathbf{d},\mathbf{d})=-2$.  It  defines a reflection $x\mapsto x+ Q(x,\mathbf{d})\mathbf{d}$ which is an isometry of $S$ with fixed hyperplane $H_{\mathbf{d}}$. All these reflections generate the \emph{Weyl group} $W_S$.

The complement inside  $C$  of all  hyperplanes $H_{\mathbf{d}}$ forms a  disjoint union of fundamental domains $D$ for the action of the reflection group $W(S)$ generated by these hyperplanes.  
For some subset $R_D$ of roots the fundamental domain $D$ can be written as  
\begin{equation} \label{eq:chamber}
D=\sett{x\in C}{Q(x,\mathbf{d})>0 \text{  for all  }\mathbf{d}\in R_D.}
\end{equation}
The  subset $R_D$  turns out to be  a set of \emph{positive roots}:  all roots   are   either positive of negative integral linear combinations of roots from $R_D$. %
Conversely any  set $R$ of positive roots determines a unique chamber $D$ for which $R_D=R$.
Choosing a different cone gives a different system of positive roots and the isometries preserving this one leads to a conjugate group. Indeed, introducing the \emph{reduced orthochronous Lorentz group}  
\[
 \ogrup{S}=\sett{g\in \ogr{S}}{g(D)= D }
\]
we have
\begin{equation}\label{eqn:latisom}
\ogrp{S}= \ogrup{S}\ltimes W(S).
\end{equation}
We also introduce the following group
\begin{equation}\label{eqn:extgrp} 
\extgrup{S}:= \sett{(\epsilon, \mathbf{g}) \in \mu_2\times \ogrup{S}}{ \overline{\mathbf{g}}=\epsilon \in \disc{S}}.
\end{equation}
\begin{rmk} \label{2group} The natural homomorphism $\extgrup{S} \to \ogrup{S}$ is   injective unless $\disc S$ is trivial or a $2$-group.
\end{rmk}

 \section{Indefinite integral binary quadratic forms}\label{sec2}

 \subsection{Pell's equation}
 Let $d$ be any positive integer.  The two Pell equations associated to $d$ are:
 \begin{eqnarray} x^2- dy^2 & =&  4:  \quad\quad\mbox{the positive Pell equation;} \label{eqn:pell1}
 \\ 
 x^2- dy^2 & =& -4:  \quad \mbox{  the negative  Pell equation.}  \label{eqn:pell2}
\end{eqnarray}
Solutions  of  \eqref{eqn:pell1}  always exist, but this is not true for \eqref{eqn:pell2}.        

There are also the  \emph{reduced} Pell equations where the right hand side has been replaced by $\pm 1$. Two cases are important for us:
\begin{enumerate}
\item $d\equiv 0 \bmod 4$. In this case $x$ is even, say $x=2\bar x$ and $(\bar x, y)$ is a solution of the reduced equation for $\frac 14 d$ if and only if $(x,y)$ is a solution of the (full) Pell equation for $d$; 
\item $d\equiv 1 \bmod 4$. In this case $x$ and $y$ have the same parity and if the{y} are both even, say $x=2\bar x$, $y=2\bar y$, then $(\bar x,\bar y)$ is a solution of the reduced Pell equation for $d$ if and only if $(x,y)$ is a solution of the full Pell equation for $d$.
\end{enumerate}
Note that if $u^2 -dv^2=1$ gives a minimal positive solution,  to any prime divisor $p$ of $d$ we can associate the unique number $\epsilon(p) = \pm 1$ such that $u \equiv \epsilon(p) \bmod p$. 
 In Table~\ref{sols}) we have collected   the smallest positive solutions $(x,y)$ of the  two positive Pell equations for some values of $d$. In the first column we put $d$ factored into primes, in the second column the minimal positive  solution $(u,v)$ is exhibited for $x^2-dy^2=1$, while in the third column the minimal  solution for $x^2-dy^2=4$ can be found.  In  the last column the numbers $\epsilon(p) $  are gathered in a  vector with  entries according to the prime decomposition of $d$. The table has been composed using \cite{Pell}.

\begin{table}[h]
\begin{center}
\begin{tabular}{|c|c|c|c|}
\hline $d$ & $N=1$ & $N=4$ & $\epsilon(p)$ \\
\hline $3$ & $(2,1$)& (4,2)  & $ {-1}$ \\
\hline $5$ & $(9,4)$& $(3,1)$ & $-1$ \\
\hline $6=2 \cdot 3 $ & $(5,2$)& (10,4)  & $(1,-1)$ \\
\hline $7$ & $(8,3)$& (16,6) & 1 \\
\hline $8=2\cdot 2 \cdot 2$ & $(3,1)$& $(6,2)$ & 1 \\
\hline $11$ & $(10,3)$& (20,6) & $-1$ \\
\hline $13$ & $(649,180)$&$(11,3)$ & $-1$ \\
\hline $17$ & $(33,8)$& $(66,33)$ & $-1$ \\
\hline $15=3\cdot 5$ & $(4,1)$ &  (8,2)   & $(1,-1)$ \\
\hline $20=5 \cdot 2 \cdot 2$ & $(9,2)$& (18,3) & $(1,-1)$ \\
\hline  $21=3\cdot7$ & $(55,21)$ & $(5,1)$ & $(1,-1)$ \\
\hline   $33=3\cdot11$ &  $(23,4)$ & (46,8)  &$(-1,1)$  \\
\hline  $35=5\cdot7$ &  $(6,1)$ &  (12,2)  &$(1,-1)$  \\
\hline   $39=3\cdot13$ &$(25,4)$   & (50,8) &$(1,-1)$ \\
\hline $44= 11 \cdot 2 \cdot 2$ & $(199,30)$& (398,60) & (1,1) \\
\hline $51=3\cdot17$ &    $(50,7)$&  (100,14)& $(-1,-1)$\\
\hline $55=5\cdot11$ &    $(89,12)$ &(178,24) & $(-1,1)$\\
\hline $104=2 \cdot 2 \cdot 2 \cdot 13$ & $(51,5)$&$(102,10)$ & (1,-1) \\
\hline $105=3\cdot5\cdot 7$ &   $(41,4)$ & (82,8)& $(-1,1,-1)$\\
\hline $165=3\cdot5\cdot 11$&   $(1079,84)$&  (13,1) & $(-1,-1,1)$\\
\hline
\end{tabular}\caption{
Minimal positive solutions of $x^2-dy^2=N$.}\label{sols}
\end{center}
\end{table}
The solutions of the positive Pell equation behave  fundamentally differently according to when $d$ is a  square or not:
 \begin{lemma} \label{pellstruct}  If $d$ is a square, the only solutions of   \eqref{eqn:pell1} are $u=\pm 2$ and $v=0$. 
  \newline
If $d$ is not a  square,  let $(U,V)$ be the  \emph{smallest} positive solution of \eqref{eqn:pell1}, i.e. $U,V>0$ are as small as possible and write  $\epsilon = \frac{1}{2}(U + V \sqrt{d})$. Then all solutions are generated by powers of $\epsilon$ in the sense that writing $\epsilon^n= \frac12(u+v\sqrt{d})$, $(u,v)$ is a new solution and all solutions can be obtained that way. \newline 
If $(U',V')$ is a minimal positive solution for \eqref{eqn:pell2} and $\eta= \frac 1 2 (U'+V'\sqrt{d})$, then $\eta^2=\epsilon$ gives the minimal solution for \eqref{eqn:pell1}, the even  powers of $\eta$ thus provide all solutions of the positive Pell equation, while the odd powers  yield all solutions to the negative Pell equation.\end{lemma} 
Note that  $\Q(\sqrt{d})$ is a quadratic extension of $\Q$ and the units of norm $1$ in  the ring of integers $\OO(\sqrt{d})$ in this field are the elements $\alpha= \frac{1}{2}(u+v\sqrt{d})$ where $(u,v)$ is  an integer solution of \eqref{eqn:pell1}. The trivial solution $(2,0)$ corresponds to $1$ in this ring. The elements of norm $-1$ correspond to the solutions of \eqref{eqn:pell2}. We shall need the following auxiliary result:
\begin{lemma} \label{convexsums} Let   Let $\eta= \frac 1 2 (U +  V\sqrt{d})$ be the unit corresponding to the minimal positive solution of  \eqref{eqn:pell1}. Then  $\eta^k = a_k \eta - b_k$, $-\eta^{-k}= c_k\eta -d_k$ with $a_k,b_k,c_k,d_k>0$ for $k>0$. In other words any  solution  of   \eqref{eqn:pell1} is expressible as  an integral linear combination of  the two solutions $\eta$ and  $-1$ with either negative or positive coefficients.
\end{lemma}
\proof  Note that $U\ge 2$ and that $\eta^2=U\eta -1$ and hence $a_2=U$, $b_2=1$. Then we have the recursive formulas $a_{k+1}=a_kU-b_k$, $b_{k+1}=a_k$. These inductively imply that for $k>0$ one has $a_k/b_k\ge 1 $. One needs to show that $a_k >0$ and $b_k>0$. The recursive formulas show that by induction we have $\displaystyle a_{k+1}= \big(\frac{a_k}{b_k} U -1\big) b_k \ge (U-1)b_k >0$ and $b_{k+1}=a_k>0$.  A similar argument applies to   $c_k$ and $b_k$.
\qed\endproof

\subsection{Isometries}
In this section we summarize the classical theory of binary quadratic forms (over a field $k$ of characteristic $0$) in a way adapted to our needs.  We are in particular interested in the associated orthogonal groups.

Recalling our convention, a quadratic form $q(x,y)=ax^2+ bxy+cy^2$  is associated to the bilinear form $Q$ given by
\[
Q(\vc{v},\vc{w})= \Tr \vc{v} \begin{pmatrix} a & b/2 \\
						b/2 & c\end{pmatrix} \vc{w}.
						\]
\begin{dfn} \label{quadrdefs}    We say that $Q$ is \emph{equivalent} to $Q'$, written $Q'\sim Q$, if  for some invertible $2$ by $2$ matrix  $P$ one has  $Q'= \Tr P Q P$. The associated quadratic forms are also said to be equivalent: $q'\sim q$. The special case $Q=Q'$  gives the  \emph{automorphs} $T$ of $q$.   If  $\det P>0$ we speak of \emph{proper equivalence}. \\
The quantity $b^2-4ac$ is the  \emph{discriminant}  $d(Q)$ of $Q$ (or of $q$) and remains invariant under equivalence.  \\
Recall  that $q$ is called  \emph{integral} if $a,b,c\in \Z$; such a form  is \emph{primitive} if  $(a,b,c)$ have no common divisors.  If $q$ admits automorphs $P$ with $\det P=-1$ we say that $q$ is \emph{ambiguous}.
\end{dfn}
 Positive discriminant  means that $Q$ is indefinite. Indeed, one has:
\begin{lemma} \label{diagonal}  In  the field $k$ the quadratic form $q$ is equivalent to the diagonal form %
$d_{4,-d}:= 4x^2-dy^2$:
\begin{equation} \label{eqn:hyp}
16a Q = \Tr   P 
		\begin{pmatrix} 4 & 0 \\
					 0  & -d   \end{pmatrix}  P, \quad 
					 P = \begin{pmatrix}  2a & b \\ 0 & 2 \end{pmatrix}                                                         
\end{equation}
\end{lemma}
Assume  that   $q$ is integral.   By  \eqref{eqn:hyp}  one can find $\ogr{Q}$ by comparing $d_{4,-d}$ and $16aQ$.  We may furthermore assume that $q$ is \emph{primitive}.  Automorphs of  the diagonal form $d_{4,-d}$ immediately lead to   solutions of  the first of the two Pell equations which are associated to $d$; indeed, using \eqref{eqn:hyp} and Lemma~\ref{pellstruct} one finds:
\begin{prop}[\protect{\cite[Th.~50,Th.~51c]{J}, \cite[Theorem 87]{D} }] \label{notperfectsqr} 
Suppose that the quadratic form $q$ is \emph{primitive} and that $a\not=0$. The   group of special  isometries   $\so Q$ of $Q$ is isomorphic to the direct product of the cyclic group $\Z/2\Z$ generated  by $-\id$ and the infinite cyclic group generated by
\begin{equation}\label{eqn:automorphs} 
\mathbf{u}:= \begin{pmatrix}\frac 1 2 [U-bV]& -cV\\
aV&  \frac 1 2 [U+bV] \end{pmatrix}.
\end{equation}
where as before, $(U,V)$  is a minimal positive solution of  the Pell equation \eqref{eqn:pell1}.
A general proper automorph   $\pm \mathbf{u}^k$ of $Q$  is (up to sign) of the form \eqref{eqn:automorphs} with $(U,V)$ replaced with  
a suitable  solution $(u,v)$ of the Pell equation \eqref{eqn:pell1} (see Lemma~\ref{pellstruct}).
\\
 If $d$ is a square, $\so{Q}=\pm\id$. \end{prop}
\begin{rmk} \label{ogrp} To see the connection with the units in the quadratic field $K=\Q(\sqrt{d})$ we proceed as follows. Let $\set{\omega_+,\omega_-}$ be the two roots of
\begin{equation}  \label{eqn:roots}
a\omega ^2+ b\omega+c  = 0,\quad \omega_{\pm}= \frac{-b \pm \sqrt{d}}{2a}.
\end{equation}
Then  in $K$ the form $q$ is equivalent to the standard hyperbolic form $h(x,y)=2xy $:
\begin{equation*} 
  8a  Q = \Tr P \begin{pmatrix} 0 & 1 \\
 1  & 0  \end{pmatrix}  P, \quad P= 2a \begin{pmatrix}   1  &  -\omega_- \\ 
                                                                                  1 & -\omega_+ \end{pmatrix}.
\end{equation*}
Moreover, we have
\[
\mathbf{u}=P^{-1} \begin{pmatrix} \epsilon & 0 \\ 0 & \epsilon^{-1} \end{pmatrix} P.
\]
This representation shows that  $\sogrp {Q}$ is generated by $\mathbf{u}$: the first  quadrant can be taken to represent the cone $C$   (see \eqref{eqn:cone})  and $\mathbf{u}$ preserves $C$ while $-\id$ does not.
\end{rmk}

To find all isometries, one needs to add at most one involution: 
\begin{prop}[\protect{\cite[Th.~52]{J}}]  \label{indirect} If $q$ is ambiguous, i.e. $Q$ admits an isometry of determinant $-1$, then  the form  is   properly equivalent to either one of the following two classes of special forms
\begin{eqnarray}     
d_{[a,c]}  & = &ax^2+c y^2 \label{eqn:diagform} \\
{\tilde d}_{[a,c]}  &=& a x^2+axy+cy^2.  \label{eqn:almostdiagform} 
\end{eqnarray}
Such forms are indeed ambiguous:  in the diagonal case put $w=0$ and in the non-diagonal  case put $w=1$, then the  involution $\mathbf{a}'=\begin{pmatrix} 1 & w\\ 0 & -1 \end{pmatrix}$  is an isometry of the corresponding bilinear form.  If  $\mathbf{a}$ is the corresponding involution for $q$, we have $\ogrp{Q}=\langle \mathbf{a},\mathbf{u}\rangle$. 
\end{prop}
\begin{rmk} \label{indirects} 
In fact, Jones does not prove that $\mathbf{a}$ preserves the components of the light cone. In the diagonal case this  can be seen as follows. 
One can take $C$ to be (smallest)  sector of the plane bounded by the two lines  $x=\pm \sqrt{-c/a} {y}$ and $\mathbf{a}'$ preserves this sector.    The non-diagonal case is similar.
\end{rmk}

Proper equivalence makes use of  \emph{adjacency}: 
we say  $q'$ is right adjacent  to $q$ (or $q$ left adjacent to $q'$) if $q'= \Tr P q P$ with $P=P_e:=\begin{pmatrix} 0 & -1 \\ 1 & e \end{pmatrix}$, $e\in\Z$. It replaces the coefficients $(a,b,c)$ of $q$ by $(c, -b+2e c , a-eb +ce^2)$. 
\begin{exmple} \label{adjacent} Using $P_0$ the form  $x^2 + b  xy + c y^2$   is adjacent to  $cx^2 - bxy + y^2$ which in turn, using $P_{-e}$  is adjacent to $x^2 + (b-2e)xy + (c- b e +e^2)y^2$. So, if $b=2e$ is even,  our  form is equivalent to the diagonal form $x^2 + (c-e^2)y^2$.  On the other hand, if $b=2e+1$ is odd, the discriminant $d$ is equivalent to $1$ modulo $4$ and the form  is equivalent to  $x^2+xy+(c-e-e^2)y^2$. In both cases, according to Prop.~\ref{indirect} there is an involution of determinant $-1$. Explicitly,
\[
\mathbf{a} = \begin{pmatrix}  1   &  b  \\
0  &  -1   \end{pmatrix}
\] 
\end{exmple}

\subsection{Discriminant groups 
of ambiguous forms
} \label{ssec:ambi}

{
Up to now we assumed that $q=ax^2+bxy+cy^2$ is primitive. For the computation of the groups $\extgrup{Q}$ we can no longer assume this. Since $Q$ is even, we can always write $Q=2nq$ with $n\in \Z$ and $q$ primitive. The discriminant form $\disc{Q}$ then is generated by the columns of the matrix
\[
  \frac{ 1}{  nd}   \begin{pmatrix} 2c & -b \\
-b & 2a
\end{pmatrix},\quad d=d(q)= b^2-4ac.
\]
This enables us to find explicit sufficient conditions for $\overline{\mathbf{u}}^k=\pm 1$ to hold. Indeed, let $(U,V)$ be the minimal positive  solution of the Pell equation~\eqref{eqn:pell1}  and let $(u_k,v_k)$ be the solution obtained by writing $\half (u_k +\sqrt{d} v_k) = \left[\half (U+\sqrt{d} V)\right]^k$. Then we have:
\begin{lemma}\label{upows}
 Suppose that the following conditions hold
 \begin{eqnarray}
 (u_{{k}}-2) c & \equiv &  0 \bmod (nd) \label{eqn:disc1}\\
 \half d v_{{k}} + (-\half u_{ {k}} +1)b   &\equiv&  0 \bmod(nd), \label{eqn:disc2}\\
 { -\half d v_k + (-\half u_k +1)b}   &  {\equiv}&   {0 \bmod(nd)} \label{eqn:disc3} \\
 { (u_k-2) a} & { \equiv} &     {0 \bmod (nd)} \label{eqn:disc4} 
 \end{eqnarray}
then $\overline{\mathbf{u}}^k= 1$. For  $n=\pm 1$ these conditions  are   always verified  with  $k=2$.
\end{lemma}
\proof 
That  \eqref{eqn:disc1}--\eqref{eqn:disc4} imply $\overline{\mathbf{u}}^k= 1$  is a direct calculation   using Prop.~\ref{notperfectsqr}. That these condition are satisfied for $n=\pm 1$  uses that $u_2= 2+d V^2$ and $v_2=UV$. Indeed  \eqref{eqn:disc1}  { 
and \eqref{eqn:disc3}
} then are immediate. %
For \eqref{eqn:disc2}, note that the left hand side equals  $\half (U- bV) dV$ and hence is divisible by $d$.  {
Equation \eqref{eqn:disc4}  can be replaced by $dv_2\equiv 0 \bmod d $ in this case which is trivially true.}
\qed\endproof
Let us now pass to ambiguous forms. }
We first consider the diagonal case $q=ax^2+cy^2$ with $a$ and $c$ coprime. Then $d(q)= -4ac$. The associated forms are
$
Q= n \begin{pmatrix} 2a & 0 \\ 0 & 2c
\end{pmatrix} 
$ with
 \[
\disc{Q}= \langle  \displaystyle \begin{pmatrix}\frac{1}{2na}\\ 0 \end{pmatrix}, \begin{pmatrix} 0\\ \frac{1}{2nc} \end{pmatrix}\rangle.
\] 
Let $(S,V)$ be the minimal positive solution of $x^2+acy^2=1,\,$ that is $(2S,V)$ is the minimal positive solution of $x^2-dy^2=4$ as in Prop.\ref{notperfectsqr}.  In this case  $ac<0$ and we may assume that $a>0, c<0 $ and $a< |c|$. We allow $n$ to be negative. 
\begin{lemma} \label{diagdisc} 
If $|n|\ge 2$ no isometry of determinant $-1$ induces $\pm \id$ on $\disc{Q}$. 
\newline
For $n=\pm  1$ and $a\not=  1$ the necessary and sufficient condition for such an involution to exist is that $V$ be even,  $S\equiv \pm 1 \bmod 2a$ and $S\equiv \mp 1 \bmod 2c$ and then always  $\overline{\mathbf{a}}'\overline{\mathbf{u}}=\pm \id$. 
\newline
For $n=\pm  1$ and $a= 1$ one has $\overline{\mathbf{a}}'=-\id$.
 \end{lemma}
\proof  We let $(s_k,v_k)$ be the solution of $x^2+acy^2=1,\,$ obtained from $(S+\sqrt{-ac}V)^k$ (cf. Lemma~\ref{pellstruct}). 
{For simplicity we write $(s,v)$ instead of $(s_k,v_k)$.  
We calculate
\begin{eqnarray*} 
\mathbf{b}&:= &\mathbf{ a' u}^k = \begin{pmatrix} s &  -cv \\
-av   & -s\end{pmatrix},\\
 \mathbf{b}\begin{pmatrix}\frac{1}{2na}\\ 0 \end{pmatrix}&=&  \begin{pmatrix}\frac{s}{2na}\\ \frac{-v}{2n} \end{pmatrix}, \quad\quad \mathbf{b}\begin{pmatrix}0 \\ \frac{1}{2nc}  \end{pmatrix}=\begin{pmatrix}\frac{-v}{2n}\\ \frac{-s}{2nc} \end{pmatrix}.
\end{eqnarray*}
\newline
Secondly, }
$\bar{\mathbf{b}}=\pm \id$ only if $n=
{\pm 1}$ and   then $v$ has to be even while $s\equiv \pm 1 \bmod 2a$, $s\equiv \mp 1 \bmod 2c$ is a necessary and sufficient condition. 
Since in this case by Lemma~\ref{upows} $\overline{\mathbf{u}}^2=1$, we can reduce to the case $k=0$ or $k=1$.
If $k=0$ we see that $a=1$ and if $a\not=1$ the necessary and sufficient conditions are as stated.
\qed\endproof
\begin{rmq} 
If   $n=\pm 1$ and $S\equiv \pm 1 \bmod (2a)$ and  $S \equiv \mp 1 \bmod (2c)$ hold simultaneously, two cases occur: if $V$ is even, then  $\overline{\mathbf{a}}'\overline{\mathbf{u}}=\pm \id$ 
but if $V$ is odd    $\overline{\mathbf{a}}'\overline{\mathbf{u}}^{\ell}$ can   never be equal  $\pm1$.

For  $a=1$ and $n=\pm 1$ the situation is different: if   $V$ is  even and $S\equiv \mp 1 \bmod (2c)$ we have $\overline{\mathbf{u}}=\pm \id$,  and if $V$ is odd   then $\overline{\mathbf{u}}^{2}= \id$.    
 \end{rmq}
 
\begin{exmples} Using Table~\ref{sols} it is easy to find examples where the conditions hold  and where these  fail: they hold for  $(a,c)= (1,-5) $,  $(1,-13)$, $(1,-17)$, 
$(1,-3)$,  $(1,-7)$, $(1,-11)$, $(3,-11)$, $(3,-13$), $(4,-5)$ but   fail  for $(a,c)=  (3,-5), (3,-7)$, $(3,-17)$, $(2,-4)$.
See  table~\ref{pows} for  more complete information. 
\begin{table}[h]
\begin{center}
\begin{tabular}{|c|c|c|}
\hline 
$(a,c)$ & smallest $k$ with $\overline{\mathbf{u}}^k=\pm\id$& ditto for $\overline{\mathbf{a}}'\overline{\mathbf{u}}^k=\pm\id$  \\
\hline $(1,-5)$ & $\overline{\mathbf{u}}=\id$& $\overline{\mathbf{a}}'=-\id$  \\
\hline $(1,-13)$ &$\overline{\mathbf{u}}=\id$& $\overline{\mathbf{a}}'=-\id$ \\
\hline $(1,-17)$ & $\overline{\mathbf{u}}=\id$& $\overline{\mathbf{a}}'=-\id$ \\
\hline $(1,-3)$ & $\overline{\mathbf{u}}^2=\id$& $\overline{\mathbf{a}}'=-\id $ \\
\hline $(1,-7)$ &$\overline{\mathbf{u}}^2=\id$& $\overline{\mathbf{a}}'=-\id  $  \\
\hline $(1,-11)$ & $\overline{\mathbf{u}}^2=\id$& $\overline{\mathbf{a}}'=-\id  $  \\
\hline $(3,-11)$ & $\overline{\mathbf{u}}=\id$& $\bar{\mathbf{a}}'\bar{\mathbf{u}}=-\id $ \\
\hline  $(3,-13)$ & $\overline{\mathbf{u}}=\id$& $\bar{\mathbf{a}}'\bar{\mathbf{u}}=-\id$ \\
\hline  $(3,-5)$ &  $\overline{\mathbf{u}}^2=\id$& none   \\
\hline   $(3,-7)$ &  $\overline{\mathbf{u}}^2=\id$& none  \\
\hline   $(3,-17)$ &$\overline{\mathbf{u}}^{
{2}}=\id$& none \\
\hline
\end{tabular}
\caption{Minimal  numbers $k$ with $\overline{\mathbf{u}}^k=\id$ and $\overline{\mathbf{a}}' \overline {\mathbf{u}}^k=-\id$.}\label{pows}
\end{center}
\end{table}
\end{exmples}
We now pass to the non-diagonal case $q=ax^2+axy+cy^2$, $(a,c)=1$. Here $Q=n \begin{pmatrix}  2a & a \\a  &2c \end{pmatrix}$, $d(q)=a(a-4c)$. Since $d(q)>0$ we may suppose $a>0$ and $a>4c$. Again, we allow $n$ to be negative. We have  
\[
\disc{Q}=\langle \begin{pmatrix}\frac {1}{na}\\0\end{pmatrix}, \begin{pmatrix}\frac{1}{n(a-4c)}\\\frac{-2}{n(a-4c)}\end{pmatrix} \rangle
\]
 and in this case:
\begin{lemma}  \label{nondiagdisc}  
If $|n| \ge 2$ and $(n,a)\not=(\pm 2,1)$ no isometry of determinant $-1$ induces $\pm \id$ on $\disc{Q}$. \newline
For $n=\pm 1, a\not= 1, 2$ the necessary and  sufficient condition for such an involution to exist is that  $U\equiv \pm 2 \bmod a$ and $U\equiv \mp 2 \bmod (a-4c)$.  If this is the case  $\overline{\mathbf{a}}'\overline{\mathbf{u}}=\pm \id$. \newline
For $(n,a)=(\pm 1, 1),(\pm 1, 2)$ or $(\pm 2, 1)$  we always have  $\overline{\mathbf{a}}'=-\id$.
\end{lemma}
\proof
For simplicity we write $(u,v)$ instead of $(u_k,v_k)$.  We find
 \begin{eqnarray*}
\mathbf{b}&:= &\mathbf{ a 'u}^k = \begin{pmatrix} \frac 12 (u+av) &  \frac 12 (u+av) -cv \\
-av   & -\frac 12 (u+av)\end{pmatrix}
 \end{eqnarray*}
 so that 
  \begin{eqnarray*}
   \mathbf{b}\begin{pmatrix}\frac{1}{na}\\ 0 \end{pmatrix} &=& \begin{pmatrix}\frac{\frac12(u+av)}{na}
   \\      {\frac{-v}{n}} \end{pmatrix}= \pm \begin{pmatrix} \frac{1}{na} \\ 0 \end{pmatrix},\\
 \mathbf{b}\begin{pmatrix}\frac{1}{n(a-4c)}\\
 \frac{-2}{n(a-4c)} \end{pmatrix}& =&  \begin{pmatrix}\frac{-\frac 12 (u+ 
 {a}v) +2cv}{n(a-4c)}\\ \frac{u}{n(a-4c)} \end{pmatrix}=
 \pm \begin{pmatrix} \frac{1}{n(a-4c} )\\ \frac{-2}{n(a-4c)} \end{pmatrix}.
\end{eqnarray*}
Now the proof proceeds as in the diagonal case.
\qed\endproof
\begin{exmples} Using Table~\ref{sols} one finds the examples gathered in Table~\ref{powsbis}.
\begin{table}[h]
\begin{center}
\begin{tabular}{|c|c|c|}
\hline  
$(a,c)$ & smallest $k$ with $\overline{\mathbf{u}}^k=\pm\id$& ditto for $\overline{\mathbf{a}}'\overline{\mathbf{u}}^k=\pm\id$  \\
\hline $(1,-1)$ & $\overline{\mathbf{u}}^2=\id$& $\overline{\mathbf{a}}'=-\id$  \\
\hline $(1,-3)$ &$\overline{\mathbf{u}}^2=\id$& $\overline{\mathbf{a}}'=-\id$ \\
\hline $(1,-4)$ & $\overline{\mathbf{u}}^2=\id$& $\overline{\mathbf{a}}'=-\id$ \\
\hline $(3,-1)$ & $\overline{\mathbf{u}}=\id$& none \\
\hline $(3,-2)$ &$\overline{\mathbf{u}}^2=\id$& none  \\
\hline $(7,1)$ & $\overline{\mathbf{u}}^2=\id$& none  \\
\hline $(21,4)$ & $\overline{\mathbf{u}}^{
{2}}=\id$& $\bar{\mathbf{a}}'\bar{\mathbf{u}}=-\id$ \\
\hline  $(15,1)$ & $\overline{\mathbf{u}}^{
{2}}=\id$&none \\
\hline  $(35,8)$ &  $\overline{\mathbf{u}}^{
{2}}=\id$&   $\overline{\mathbf{a}}'\bar{\mathbf{u}}=-\id$   \\
\hline
\end{tabular}
\caption{Minimal  numbers $k$ with $\overline{\mathbf{u}}^k=\id$ and $\overline{\mathbf{a}}' \overline {\mathbf{u}}^k=-\id$.}\label{powsbis}
\end{center}
\end{table}
\end{exmples}

 \subsection{Roots} \label{ssec:roots}
There is a general theory of  representations of integers by $q$. See \cite[\S~46]{D}. We only need the theory of representations of $-1$. These correspond to representations of $-2$ by $2Q$, i.e. to the roots of the corresponding even lattice.
  This general prescription yields:
 \begin{lemma}  \label{reprone} Suppose that $(x,y)$ is a solution  for $q(x,y)=-1$ with relative prime integers $x$ and $y$. Then either  
\\ 
\noindent {\rm(*)} $d\equiv 0 \bmod 4$ and  $q\sim d_{[ -1, \frac 14 d] }$  hence, if $(u,v)$ is a positive solution for \eqref{eqn:pell1}, then $(\frac 1 2 u,v)$ gives a positive representation of  $-1$ for $d_{[ -1, \frac 14 d] }$ and conversely.  
 \newline
 or 
 \\
 \noindent {\rm (**)} $d\equiv 1 \bmod 4$ and $q\sim {\tilde d}_{[ -1, \frac 14 (d-1)] }$ and $(\frac{1}{2}(u-v), v))$ gives  a positive representation for  $-1$ for ${\tilde d}_{[-1, \frac 14 (d-1)]}$ and conversely.
  \end{lemma} 
  \begin{rmk}\label{rootsambig} By Prop.~\ref{indirect} this Lemma implies that the quadratic forms $q$ for which $2Q$ has roots are automatically ambiguous.
  \end{rmk}
  
 We study these two forms in more detail.
   \begin{lemma}   \label{repminone} Let $(U,V)$ the minimal positive solution for \eqref{eqn:pell1}.  The   solutions of  the equation  $d_{[-1, \frac 14 d]} (x,y)=-1$ are either positive or negative linear combinations of the  two \emph{basic solutions}  $\mathbf{e}=(-1,0)$ and $\mathbf{f}=(\frac 1 2 U, V)$.
 The solutions of the equation $\tilde d_{[-1, \frac 1 4 (d-1)]}=-1$  are either positive or negative linear combinations of the  two \emph{basic solutions}  $\mathbf{e}=(-1,0)$ and $\mathbf{f}=(\frac 1 2 (U-V), V)$.
In both cases the involution $-\mathbf{a'u}$  generates the group $\ogrup {2Q}$
\newline
{The group  $\extgrup{2Q}$  is trivial unless  \eqref{eqn:pell2} is solvable and then the involution $(-\id,-\mathbf{a'u})$ generates $\extgrup{2Q}$.}
 \end{lemma}
\proof
The first assertion follows from Lemma~\ref{convexsums} since in both cases $1\in \OO(\sqrt{d})$  corresponds to $(1,0)$, while  $\eta \in \OO(\sqrt{d})$ corresponds to
$(\frac 1 2 U, V)$ in the first case and to  $( \frac 1 2 (U-V),V)$ in the second case.
\newline
We only treat the diagonal case $q=d_{[-1, \frac 14 d]}$; in the non-diagonal case the computations are similar.   Recall the notion of positive root from \S~\ref{sec1}. Here we can take  the two basic roots  as   positive roots with respect to $2Q$.  Since $\ogr{2Q}$ consists of elements of the form   $\pm \mathbf{u}^k$ or $\pm \mathbf{a'} \mathbf{u}^k$  it suffices to determine which of    these elements  preserve the set of  basic roots. By direct computation one finds that $\mathbf{a'u}(\mathbf{e})=-\mathbf{f}$ and $\mathbf{a'u}(\mathbf{f})=-\mathbf{e}$ and hence  $-\mathbf{a'u}$ preserves the cone $D$ defined in \eqref{eq:chamber}. By a similar computation one sees that all the other elements of $\ogr{2Q}$ do not preserve the set of basic roots. The first statement  follows.
The second assertion follows from  Lemma~\ref{pellstruct}. Indeed, since the negative Pell equation has a minimal solution  of the form $(2s,t)$, the minimal solution for the positive Pell equation is of the form $(U,V)=(4s^2+2,2st)$  and since $V$ is even and $U\equiv -2  \bmod d$ one calculates directly, as in the proof of Lemma~\ref{diagdisc}  that  $\bar{\mathbf{a}}'\bar{\mathbf{u}}=\id$. Conversely, suppose that  $\mathbf{a'u}$ induces $ \pm \id$.  Again, as in the proof of Lemma~\ref{diagdisc}, we  must have  that $V$ is even and $U\equiv \mp 2 \bmod d$. In this case write $\frac 12 U= \mp 1+   \frac 12 k d, V= 2 m$. Then $m^2=k( \frac 14 dk \mp  1)$ shows that $k=t^2$ and $\frac 14 dk \mp  1=s^2$ are squares of integers $(s,t)$  for which then $s^2 - \frac 14 dt^2=\mp 1$  {(remember that $d$ is divisible by $4$ in this situation so that  $k$ and $\frac 14 dk \mp  1$ are coprime)}. Note that  the plus sign is excluded since $(U,V)$ was supposed to be a minimal solution.  So $\mathbf{a'u}$ induces $ \id$ and $(s,t)$ solves the negative Pell equation.
\qed \endproof 
  The drawback of Lemma~\ref{reprone} is that it is hard to apply in general. But Lemma~\ref{repminone}  suggest a further link with the solvability of  \eqref{eqn:pell2}. Indeed, if $q$ has leading coefficient $1$ one  verifies immediately:
\begin{lemma} \label{rootspecialcase} Let $q=x^2+bxy+cy^2$ with discriminant $d=b^2-4c$. Then
 \eqref{eqn:pell2} is solvable if and only if  $q$ represents $-1$.  
 Indeed, a solution $(u,v)$ of  \eqref{eqn:pell2} yields a solution
  \begin{equation*}  
  \vc{v}= 
  \begin{pmatrix}\frac{1}{ 2} [u-bv]\\
 v  \end{pmatrix} 
  \end{equation*}
of  $q(\vc v)= -1$ and conversely.
\end{lemma}
 As  in the proof of Lemma~\ref{repminone}  to test solvability of  \eqref{eqn:pell2} it can be useful to investigate the minimal  solution of the positive Pell equation. This is illustrated by the following example.    
 
 \begin{exmple} \label{exmple1}
 Consider  the form $q_\delta:=x^2+\delta xy + y^2$ with discriminant $d=(\delta^2-4)$.  Then $(U,V)=(\delta,1)$ gives the smallest positive solution of \eqref{eqn:pell1} and if a solution for the negative Pell equation would exist we would have ${U'}^2+(\delta^2-4){V'}^2=2\delta $ and $U' \cdot V'=1$ which forces $\delta=3$. Hence, unless $\delta=3$, the form $q_\delta$ does not represent $-1$.  We conclude that the associated even form with matrix
\[ 
 Q'_\delta:= \begin{pmatrix}
 2 & \delta\\
 \delta &  2\end{pmatrix}
 \]
does not  represent  $-2$ unless $\delta=3$. So $\ogrp Q=\ogrup Q$ in this case and the group is generated by $\mathbf{u}$ and $\mathbf{a}'$.
  \end{exmple}
 
\subsection{The discriminant $d$ is a   square}  \label{repszero} 
An integral form $q=ax^2+bxy+cy^2$ represents $0$ if and only if $d(q)$ is a  square:  $d(q)=\delta^2$ with $\delta\in \Z$. Indeed, this follows immediately from Lemma~\ref{diagonal}.  We suppose that $q$ is primitive. We have to consider all even multiples of  $q$.
We prove here:
\begin{prop} \label{discsquare}{\rm 1)}  If $d$ is a square and $q$ is not ambiguous, then  $\ogrup {2nq}=\extgrup {2nq}=\id$. \\
{\rm 2)}  If $q$ is ambiguous and equivalent to the diagonal form there are three cases: {\rm 2a)} $q=(x^2-y^2)$,  {\rm 2b)} $q=(x^2 -a^2y^2)$ $a\ge 2$ , {\rm 2c)} $q= (a^2x^2-y^2)$, $a\ge 2$.
In all cases $\ogrup {2nq}$ is cyclic of order $2$ while $\extgrup {2nq}=\id$ except in the case when $n=1$ and $a=1$. Then $\ogrup {2q} =\id$, but $  \extgrup { q}=\mu_2$.
\\
{\rm 3)} In the non-diagonal case we must have $q(x,y)=\pm (ux+vy)(ux+(u-v)y)$ with $u\ge 1$ and $u$ and $v$ coprime. So $q= \tilde{d}_{u^2,v(u-v)}$.
We have the following cases:  {\rm 3a)} $2q$ is unimodular:  $u=1,v=0$ or $u=v=1$ and then $2q$ is equivalent to the hyperbolic plane, {\rm 3b)} $u=1$ but $v\not=0,1$, {\rm 3c)} $u\not =1$.  In case {\rm 3a)}  $\ogrup {2nq}$ is cyclic of order two generated by the involution $ \begin{pmatrix}  0 & 1 \\1 & 0\end{pmatrix}$, but only for $n=1$or $n=2$  this involution induces $\id$ on $\disc{nq}$. In case {\rm 3b)} the group  $\ogrup {2nq}=\extgrup {2nq}$ is cyclic of order $2$ for $n=1$ or $n=2$. In the remaining cases  $\ogrup {2nq}$ is cyclic of order $2$ but $\extgrup {2nq}=\id$. 
\end{prop}
\proof
Prop.~\ref{notperfectsqr} implies the results in  case $q$ is not ambiguous, since then $\ogr q= \so q =\pm \id$. 
Next consider the case  of an ambiguous form $q$. By Prop.~\ref{indirect}  $q$   must be  equivalent  to $d_{[a,c]}$ or to  ${\tilde d}_{[a,c]}$. 
\newline 
\underline{Diagonal case}. Since $d$ is a square we have the three cases as stated. The group of isometries equals $\set{\pm \id, \pm \mathbf{a}= \pm \begin{pmatrix}  1 & 0\\ 0& -1\end{pmatrix} }$. The    isotropic vectors are  the multiples of $\mathbf{f}_{\pm}=(1,\pm 1)$ (in case a), $\mathbf{f}_{\pm}=(a,\pm 1)$ (in case b),  $\mathbf{f}_{\pm}=(\pm 1,a)$  (in case c)  and $\mathbf{a}$ interchanges $\mathbf{f}_+$ and $\mathbf{f}_-$ The cone $C$ is bounded by the  half lines     $\R_+\cdot \mathbf{f}_+$  and $\R_+\cdot \mathbf{f}_{-}$   and $\mathbf{a}$ preserves $C$.  Hence $\ogrp q= \ogrup q= \set{\id, \mathbf{a}}$ unless $2q$ has roots. The latter  is only the case if  $q=2x^2-2y^2$,  and then the roots are  $\mathbf{d}=\pm (0,1)$. In this case the line $H_{\mathbf{d}} $ orthogonal to $\mathbf{d}$ divides $C$ in two half cones  $C^\upvect$ and  $\mathbf{a} C^\upvect$ hence $\ogrp q= \ogrup q=\id$. In this case  $\disc{2q}=\Z/2\Z\oplus \Z/2\Z$ and hence $  \extgrup {2q}=\mu_2$ in this case.
If  $d\not=4$ we have that  $\disc{2nq}\not=1$ is not a $2$-group and hence   $\extgrup{2nq}=\id$ (compare with Remark~\ref{2group}).
\newline
\underline{Non-diagonal case}.   In the second case there exists $u,v \in\Z$ for which    $a= \pm  u^2 $, $c=\ \pm  v(u-v)$.  The form is equivalent to $(ux+vy)(ux+(u-v)y)$.   The two independent isotropic vectors $(v,-u)$ and $(v-u,u)$ span $C$ and are interchanged by $\mathbf{a}$. Hence  $\ogrp q= \set{\id, \mathbf{a}}$.  

As before, $2q$ only has roots  if   $q\sim x^2+xy$, or, equivalently $2q\sim H$. This is case a). For the hyperbolic plane  the roots are $\pm \mathbf{d}$ with $\mathbf{d}=(1,-1)$  and the line $H_{\mathbf{d}}$ divides $C$ in $C^\upvect$ and  $\mathbf{a} C^\upvect$ so that $\ogrp Q=\set{\id}$ in this case. However $\extgrup H= \mu_2$ also $\extgrup {2H} =\mu_2$ while $\extgrup {nH}=\id$ for $n\ge 3$.

In case b) we have $2nq=\pm  n\begin{pmatrix} 2 & 1 \\ 1 & -2v(v-1)\end{pmatrix}$ and $\disc{2nq}$ is generated by $ \displaystyle \frac{1}{n(2v-1)^2}\begin{pmatrix}  -1 \\2 \end{pmatrix}$ and 
$ \displaystyle \frac{1}{n}\begin{pmatrix} 1  \\0  \end{pmatrix} $ . One sees that $\overline{\mathbf{a}}=-\id$ if $n=1$ or $n=2$ so that then  $\extgrup{2nq}$ is generated by the involution $(-\id,\mathrm{a})$. If $n\ge 3$ we have $\extgrup{2nq}=\id$.

In case c) we easily find that $\overline{\mathbf{a}}\not= \pm\id$ so that  $\extgrup{2nq}=\id$.
 \qed\endproof

\subsection{ The discriminant $d$ is not a  square}\label{repsnotzero}

\begin{prop} \label{notasquare} If $Q$ has no roots  and $q$ is not ambiguous $ \extgrup Q$ is   infinite cyclic.  If $Q$ has no roots and $q$ is  ambiguous $ \extgrup Q$ is either infinite cyclic or the infinite dihedral group $\bD_\infty$.
If $Q$ admits roots,   $ \extgrup Q$ is cyclic of order two if  the negative Pell equation \eqref{eqn:pell2} for $d=d(q)$ is solvable  and trivial otherwise.
\end{prop}
\proof
Let us first consider the case when $q$ is not ambiguous, i.e.  $\so{Q}=\ogr{Q}$. By Prop.~\ref{notperfectsqr} this group is a product of $\pm\id$ and the cyclic group generated by  $\mathbf{u}$. In this case there are no roots, since roots are only possible for ambiguous forms (see Remark~\ref{rootsambig}). Hence by Remark~\ref{ogrp} we have that $\ogrp{Q}= \ogrup Q $ is the cyclic group with  $\mathbf{u}$ as generator.  It follows that some power of $\mathbf{u}$  acts as $\pm \id$ on the discriminant group and hence in this case $ \extgrup Q$ is an infinite  cyclic group with this generator.

Next, suppose that  $q$ is ambiguous. Then $\ogr{Q}$ is generated by $-\id$, $\mathbf{u}$ and the involution $\mathbf{a}$ corresponding to $\mathbf{a}'$ (see Prop~\ref{indirect}). One has $\mathbf{aua}= \mathbf{u}^{-1}$ as  a generating relation.   Hence  the group $\ogr{Q}$ consists of the elements $\pm \mathbf{u}^k$, and $ \pm\mathbf{a u}^\ell$, where $k,\ell\in \Z$. 

First assume that there are no roots. Then $\ogrup Q= \ogrp Q$ and using Remark~\ref{ogrp} and \ref{indirects}  we see that $\ogrp Q$  is generated by $\mathbf{u}$ and $\mathbf{a}$.  There is  a second involution $\mathbf{au}$ and the automorphism group $\ogrp{
{Q}}=\ogrup{
Q}$ is the group $\bD_\infty $ with generators $\set{\mathbf{a}, \mathbf{au}}$.  One next has to study  the effect  of $\mathbf{u}$ and $\mathbf{a}$ on the group $\disc{
Q}$.   Since this group is finite,  there will be a smallest  positive integer $k$ for which $\mathbf{u}^k$ induces $\pm\id$ on $
{\disc{Q}}$. Now either  there is a smallest  positive integer $\ell$ for which $\mathbf{a} \mathbf{u}^\ell$ induces $\pm \id$ or such an integer does not exist.  In the former case $\mathbf{a} \mathbf{u}^\ell$ and $\mathbf{a} \mathbf{u}^{k+\ell}$ generate the subgroup of isometries inducing $\pm \id$ on $\disc{S}$.  So  $\extgrup Q$ is the free product generated by two distinct  
{non} commuting involutions.  In the latter case the group   $\extgrup Q$ is cyclic generated by  $\mathbf{u}^k$.

Finally the case that $q$ has roots. Then Lemma~\ref{repminone} gives the answer.
  \qed\endproof

\section{Automorphism groups of K3-surfaces} \label{sec3}
For the results in this section   we refer to \cite{BPV} and \cite{PSS}.  From now on  the intersection product between two classes $c,d$ in the lattice  $L$  will be written $c\cdot d$ instead of $Q(c,d)$ and we write $c^2$ instead of $c\cdot c$.

Let $X$ be an algebraic K3-surface, i.e. a simply connected projective surface with trivial canonical bundle. Up to multiplicative constants there is a unique holomorphic $2$-form $\omega_X$ on $X$ and it is nowhere zero. 
The second cohomology integral group   equipped with the unimodular intersection form     is known to be isometric to the unique even unimodular lattice  $(L,Q)$ of signature $(3,19)$.  Any  choice of such an isometry (a \emph{marking}) identifies the line $\C\cdot \omega_X \subset H^2(X;\C)$ with a a line  $\C\cdot \omega\subset  L\otimes \C$, the \emph{period} of $X$.
An automorphism $g$  of $X$ induces an  automorphisms of $L$ which preserve the complex line $\C\omega$. So $g$ preserves  the lattice $S= \omega^\perp \cap L$  which corresponds to the Picard lattice $S_X$ as well as its orthogonal complement $T=S^\perp\subset L$, which corresponds to  the \emph{transcendental lattice}.

First consider the action on $S$.
It has signature $(1,r)$. Hence, the description of \S~\ref{sec1} applies to $\ogr S$. The canonical choice for the subcone $ C^\upvect $ is the ample cone which corresponds to the  effective roots:
\[
C^\upvect =\sett{x\in C^+}{ x\cdot d >0,\quad \mbox{for all effective roots } d}.
\]
The automorphisms of $X$ preserve this subcone so that $\aut X$ acts on $S$ as a subgroup of $\ogrup S$.  From the Torelli theorem one can deduce the complete description of the automorphism group as follows.
\begin{thm} Let $X$ be a K3-surface. Choose  an isometry $H^2(X)\mapright{\sim} L$.  The group $\aut X$ corresponds to the subgroup $G\subset \ogr L$ consisting of those $g\in \ogr L$ which preserve the period of $X$ and  for which  $g|S\in \ogrup S$.  
\end{thm}
To determine the full group $G$ it suffices to know  its restriction  to $S$ and $T$.  Since $L$ is unimodular, the groups $\disc S$ and $\disc T$ are naturally isomorphic. An automorphism of $L$ induces the same automorphism on both groups. Conversely, a pair $(g_1,g_2)\in \ogr{S}\times\ogr{T}$ can be lifted to an automorphism of $L$ if $g_1$ and $g_2$ induce the same automorphism on  $\disc  S \simeq \disc T$. In particular, if $g_2=\pm \id$, this states that  any $g_1\in \ogr{S}$  which induces $\pm \id$  on $\disc S$ lifts to a unique isometry of $L$ restricting to  $\pm \id$  on  $T$.

 The first condition for $g\in\ogr L$ to belong to $G$ reads   $g(\omega)= \lambda \omega$ for some root of unity $\lambda$.  In particular, the  automorphisms of $X$ acts as a finite group on $T$.
 If $\lambda\not=\pm 1$  the period point  $\omega$ is an eigenvector of a non-trivial isometry of $L$. This imposes algebraic conditions on $\omega$ and hence  for very general $\omega$ such automorphisms cannot exist.  This leads to:
 \begin{dfn}  An algebraic K3-surface $X$ is \emph{general with respect to automorphisms}, or Aut-\emph{general} if its automorphisms  preserve its period up to sign.
\end{dfn}
  It is equivalent to the statement that the automorphisms act as $\pm\id$ on the lattice $T$.
Hence:
\begin{prop} Let $X$ be an \emph{Aut}-general algebraic  K3-surface. Choose a marking to identify $H^2(X)$  with the lattice  $L$ and let $S$ be the sublattice of $L$ which corresponds to the Picard lattice of $X$. The automorphism group of $X$  can be identified with the group $\extgrup{S}$ (see \eqref{eqn:extgrp}).
\end{prop}
As a consequence of the previous discussion and the results in \S~\ref{sec2} we have:
\begin{corr} For $X$ a  K3-surface with Picard number $2$ the group $\aut  X$ is finite precisely when the Picard lattice  $S_X$ contains divisors $L$ with $L^2=0$ or with $L^2=-2$. If moreover, $X$ is Aut-general  we have  in this situation that $\aut  X=\id$ or $\aut X$ is cyclic of order $2$.   See Prop.~\ref{repszero} and Lemma~\ref{repminone}  for more details.

If $S_X$ does not contain such divisors and if moreover $S_X$ is not ambiguous, then  $\aut X$ is infinite cyclic, but  if $S_X$ is ambiguous,  then   $\aut X$ is either infinite cyclic or the infinite dihedral group $\bD_\infty$.
\end{corr}
\begin{rmq} This  reproves the finiteness criterion from \cite[\S~6]{PSS}.  
\end{rmq}

To extend the  result to   K3-surfaces which are not Aut-general we consider the action of  $\aut X$ on $T$. The group  $G$ acts  on $\C\cdot \omega$   through a  finite cyclic group $G'$ of order say $d$   as  a character.   On the other hand, the rational vector space $T\otimes\Q$ is a  $G'$-representation space   which  splits into a number of copies of the unique $\varphi(d)$-dimensional irreducible $\Q$-representation  space. Here $\varphi(d)$ is the Euler function. This imposes restrictions on  $G'$ .

\begin{exmple} Let rank$(S)=2$. Then  rank$(T)=20$ and hence  $\varphi(d)$ divides $20$ and as in \cite[Thm. 10.1.2]{N3} we find:

\textit{If  there is a non-trivial kernel $\ogrup L  \to \ogrup S$, 
the discriminant group $\disc S$  belongs to  the following list:  $\Z/2\Z$, $(\Z/2\Z)^2$, $(\Z/2\Z)^3$, $\Z/3\Z$, $\Z/5\Z$, $(\Z/5\Z)^2$ or $\Z/11\Z$. }

In particular, if the discriminant of $S$ is different from $2, 4, 8,3, 5,25$ or $11$ the automorphism group of any K3-surface whose Picard lattice is $S$ is the  group
$\extgrup{S}$  we just introduced and in these cases the K3-surface is automatically Aut-general.
\end{exmple}

\section{Examples}\label{sec4}

By \cite[Thm.~1.14.4]{M}    every even lattice   of 
signature   $(1,r)$ with $r\le 9$  occurs as the 
Picard lattice $S$ of some algebraic K3-surface 
and the primitive embedding $S \into L$ is unique.  In particular, all indefinite even lattices of rank $2$ may occur.   In general however it is not so easy to describe the surfaces geometrically.

We now give  some examples illustrating the theory of  \S~\ref{sec2}  some of which can be found in the existing literature.
\begin{enumerate}
\item Consider the hyperbolic lattice $H$. 
It represents $0$ as well as  $-2$.  We have seen  (Prop.~\ref{discsquare}) that $\ogr{H}$ is the identity and $\extgrup{H}=\mu_2$ which means that  the automorphism group of a general K3-surface with $H$ as Picard lattice is generated by an involution acting trivially on the Picard lattice but as $-\id$ on the transcendental lattice. The surface has a unique elliptic fibration over $\bP^1$ with a $(-2)$-section. See \cite[\S~5.4]{Ge} 
\item More generally, consider $\Lambda_{b,c}=\begin{pmatrix}  0 & b \\ b & 2c\end{pmatrix}$. This matrix represents zero and $\ogrup{\Lambda_{b,c}}$ is either trivial or cyclic of order two.  The group  $\extgrup{\Lambda_{b,c}}$ is trivial unless $(b,c)=(1,0),(2,0)$ or $(b,1)$, $b\ge 2$. This follows immediately from Prop.~\ref{discsquare}.  The corresponding K3-surfaces have been studied in \cite{Ge} where   interesting projective models are exhibited. For instance $\Lambda_{2,0}$ is realized by a double cover of $\bP^1\times \bP^1$ branched in a curve of bidegree $(4,4)$. The rulings give two elliptic pencils and the covering involution is the unique non-trivial automorphism.
 \item  Consider  $
 \begin{pmatrix}
 2 & 4\\
4 &  2\end{pmatrix}$. This matrix is equivalent to the diagonal form $d_{2,-6}= \langle2 \rangle \oplus  \langle -6\rangle$.  For this form  $\mathbf{u}=\begin{pmatrix}
 2 & 3\\
3 &  1\end{pmatrix}$ and $\mathbf{a}= \begin{pmatrix}
 1 &0\\
0 &  -1\end{pmatrix}$. These preserve $C$. The action on the discriminant group is as follows: $\mathbf{u}^2$ acts as the identity while $\mathbf{a}$ acts as $-\id$. So $\extgrup {d_{2,-6}}$ is generated by the two commuting involutions  $\mathbf{au}^2$   and $\mathbf{a}$ which both act as $-\id$ on the transcendental lattice. In fact this example has been treated geometrically. See \cite{W} where it is shown that the corresponding K3-surface $X$ is a complete intersection inside $\bP^2\times \bP^2$ of bidegree $(1,1)$ and $(2,2)$. The two  projections onto the factors $\bP^2$ realize $X$ as a double cover and the two involutions induce $\mathbf{au}^2$   and $\mathbf{a}$ on the Picard lattice.

 \item More generally $
 Q'_\delta:= \begin{pmatrix}
 2 & \delta\\
 \delta &  2\end{pmatrix}$, an  example from \cite{GL}.
 In Example~\ref{exmple1} we have seen that unless $\delta=3$ this form   represents neither $0$ nor $-2$.  If $\delta\not=3$  one has
 \[
\mathbf{u}=  \begin{pmatrix}  0 & -1 \\ 1 & \delta\end{pmatrix}, \quad \mathbf{a}=\begin{pmatrix} 1 & \delta \\0 & -1\end{pmatrix}
 \]
The group  $\disc{Q'_\delta}$ is    generated by    $\mathbf{e}=\displaystyle \frac{1}{d} \begin{pmatrix} -2  \\  \delta \end{pmatrix}\in \Z/d\Z$, or, equivalently,  by  $\mathbf{f}= \displaystyle \frac{1}{d}  \begin{pmatrix} -\delta  \\  2 \end{pmatrix} \bmod \Z/d\Z$. One verifies that $\mathbf{u}\mathbf{e}=  - \mathbf{f}$, so that $\mathbf{u}\not=\pm 1$, but $\mathbf{u}^2=1$.  On the other hand  $\mathbf{a}$ acts as $-1\in (\Z/d\Z)^\times$. It follows that the automorphism group of the K3-surface is isomorphic  to the group generated by $\mathbf{u}^2$ and $\mathbf{a}$. The first preserves the period while the second sends it to its opposite.
 \item Consider $
 Q_\delta:= \begin{pmatrix}
 2 & \delta\\ \delta &  -2\end{pmatrix}$.  In this case $(\delta,1)$ is the smallest positive solution for  the negative Pell equation and corresponds to the matrix
  \[
\mathbf{u_-}=  \begin{pmatrix}  0 & 1 \\ 1 & \delta\end{pmatrix}
\] and
\[
\mathbf{u}=   \mathbf{u_-}^2,\quad    \mathbf{a}=\begin{pmatrix} 1 & \delta \\0 & -1\end{pmatrix}
 \]
 This matrix represents $-2$: take $x=0$, $y=1$. It follows that the automorphism group of the K3-surface is finite.  In fact, by Lemma~\ref{repminone}  it is generated by $-\textbf{au}$ and hence cyclic of order $2$. This example has been treated in \cite{GL}.
 \item  The \textbf{diagonal form} $Q=2nq$ with  $q=ax^2+by^2$, $(a,b)=1$. These represent all  the \textbf{ambiguous forms with $d(q)\equiv 0 \bmod 4$}.
 \newline
 Bini considers the case $(a,b)=(d,-1)$. We find back his main result \cite[Theorem 1]{Bi}. Indeed  the form represents zero if $d$ is a square and then by Prop.~\ref{discsquare} $\extgrup{Q}$ is trivial unless $n=d=1$ in which case it is cyclic of order two. If $d$ is not a square and $n=1$ there are roots and the group $\extgrup{Q}$ is cyclic of order two by Lemma~\ref{repminone}.    If however $n\ge 2$ the group $\extgrup{Q}$ is always infinite cyclic by Prop~\ref{notasquare}. There are several projective models described in loc. cit. For example the complete intersection of bidegree $(2,1)$ and $(1,2)$ in $\bP^1\times \bP^3$ has $n=2$ and $a=1$.
\newline 
 For other values of the numbers $(a,b)$ the reader should look at Table~\ref{pows}.
 
 \item The forms $Q=2nq$, with $q$ \textbf{ambiguous and  $d(q) \equiv 1 \bmod 4$}. Such a form is equivalent to $\tilde q_{[a,c]}=ax^2+axy+cy^2$ , $(a,c)=1$: for some values of $(a,c)$ Table~\ref{powsbis} give some results.
 
 \item Forms $Q=2nq$ with $q=x^2+bxy+cy^2$ a \textbf{monic form}. By Example~\ref{adjacent} such $q$  is either equivalent to a diagonal form $q_{[1,c']}= x^2+c'y^2$, and such forms  are covered by Bini's results \cite{Bi}, or $q$ is equivalent  to $\tilde q_{[1,c']}= x^2+xy+c'y^2$. By Lemma~\ref{rootspecialcase} the latter form represents $-1$ if and only \eqref{eqn:pell2} has a solution. This gives roots if and only if $n=2$ and by Lemma~\ref{repminone} $\extgrup{Q}$ is then cyclic of order $2$. If  \eqref{eqn:pell2} has no solution Lemma~\ref{nondiagdisc} gives a recipe for determining $\extgrup{Q}$.
  \end{enumerate}

\bigskip
\noindent  	 
Dipartimento di Matematica, Universit\`a di Torino, Via Carlo Alberto n.10, 10123 Torino, ITALY

\noindent
{\it e-mail:} federica.galluzzi@unito.it

 \noindent
{\it e-mail:} giuseppe.lombardo@unito.it

\medskip
\noindent
Universit\'e de Grenoble I,
D\'epartement de Math\'ematiques
Institut Fourier, UMR 5582 du CNRS,
38402 Saint-Martin d'H\`eres Cedex, FRANCE

\noindent
{\it e-mail:} chris.peters@ujf-grenoble.fr


\begin{thebibliography}{XXXXX}

\bibitem[BHPV]{BPV} Barth, W., Hulek, K., Peters. C. and A. Van de Ven: \textsl{Compact complex surfaces}, second enlarged edition, Springer Verlag 2004.

\bibitem[Bi]{Bi}   Bini, G.: On automorphisms of some K3 surfaces with Picard Number Two, MCFA 
Annals (2005).



 \bibitem[Dickson]{D} Dickson, L. E.: \textsl{Introduction to the theory of numbers}, Dover Publ. Inc. New York 1954
 
 \bibitem[G-L]{GL}   Federica Galluzzi, F. and G.  Lombardo: On automorphisms groups of some K3 surfaces, preprint {\tt arXiv:mathAG0610972}
 
  \bibitem[G-S1]{GS}   Garbagnati, A. and A. Sarti: Symplectic automorphisms of prime order on K3 surfaces, J. of Algebra \textbf{318}  (2007) 323--350
 
  \bibitem[G-S2]{GS2}   Garbagnati, A. and A. Sarti:  Elliptic fibrations and symplectic automorphisms on K3 surfaces, preprint {\tt arXiv:matAG0801.3992}
  
  \bibitem[Ge]{Ge} Geemen, B. van: Some remarks on Brauer groups of K3 surfaces, Adv. in Math. \textbf{197} (2005) 222--247
 
 \bibitem[Jones]{J} Jones, B. W.: \textsl{The arithmetic of quadratic forms}, The Carus Mathematical Monographs \textbf{10}, MAA, John Wiley Sons, 1950.
 
 \bibitem[Kon1]{Kon1}  Kond\=o, S.: The maximum order of finite groups of automorphisms of K3 surfaces, Am. Math. J. \textbf{121} (1999), 1245--1252
 \bibitem[Kon2]{Kon2}  Kond\=o, S.: Niemeier lattices, Mathieu groups, and finite groups of symplectic automorphisms of K3 surfaces. With an appendix by Shigeru Mukai, Duke Math. J. \textbf{92} (1998) 593--603
 
 \bibitem[Mor]{M}   Morrison, D. R.: On K 3 surfaces with large Picard number,  Invent. Math. \textbf{75} (1984) 105--121 
 
 \bibitem[Muk]{Muk} Mukai, S.: Finite groups of automorphism of K3 surfaces and the Mathieu group, Inv. Math. \textbf{94} (1988) 183--221

\bibitem[Nik1]{N1} Nikulin,  V.:  Integral symmetric bilinear forms and some of their applications, Math.  USSR Izv. \textbf{14} (1980) 103--167

\bibitem[Nik2]{N2} Nikulin,  V.: Finite groups of automorphisms of K\"ahlerian $K3$ surfaces. (Russian) Trudy Moskov. Mat. Obshch. \textbf{38} (1979), 75--137

\bibitem[Nik3]{N3} Nikulin,  V.: Factor groups of automorphisms of hyperbolic forms with respect to subgroups generated by $2$-reflections. Algebrogeometric applications,  J. Soviet Math. \textbf{22}, 1401--1476, (1981)


\bibitem[P-SS]{PSS} Pjate\v ckii-Shapiro, I. and I. Shafarevi\v c: A Torelli theorem for algebraic surfaces of type K-3, Izv. Akad. Nauk. SSSR. Ser. Math. \textbf{35} (1971) 503--572.

\bibitem[Sev]{Sev} Severi, F. :Complementi alla teoria della base per la totalit\`a delle curve di una superficie algebrica, Rend. Circ. Mat. Palerma \textbf{30} (1910), 265--288



\bibitem[W]{W} Wehler, J.: $K 3$-surfaces with Picard number $2$. Arch. Math. \textbf{50} (1988), 73--82
 
\bibitem[Pell]{Pell} Quadratic diophantine equations and fundamental unit BCMATH programs,  
%
\\
\texttt{http://www.numbertheory.org/php/PELL.html}
\bibitem[Se]{Se} Serre, J.-P.: \textsl{A course in arithmetic}, Springer Verlag,
Berlin etc. (1973)
%
 \end{thebibliography}
\end{document}